\documentclass[12pt]{amsart}
\usepackage{amsthm,amsfonts, amssymb, amscd}
\newtheorem{prop}{Proposition}[section]

\newtheorem{def-prop}{Definition-Proposition}
\newtheorem{theorem}{Theorem}[section]
\newtheorem{corthm}{Corollary}[section]

\newtheorem{lemma}{Lemma}[section]
\theoremstyle{definition}

\theoremstyle{remark}

\newtheorem{example}{Example}[section]
\def\endproof{$\square$}

\def\1{1}
\def\Adj{\mbox{Adj }}

\def\1{{\mathbb I}}
\def\R{{\mathcal R}}

\def\cS{{\mathcal S}}

\def\g{{\gamma}}

\def\kbar{\overline{k}}

\DeclareMathOperator{\GL}{GL}
\DeclareMathOperator{\PGL}{PGL}
\DeclareMathOperator{\SL}{SL}

\DeclareMathOperator{\stab}{stab}
\DeclareMathOperator{\Gal}{Gal}

\begin{document}

\title[Good representations and solvable groups]
{Good representations and solvable groups}
\author{Dan Edidin}
\address{Department of Mathematics\\ University of Missouri\\
Columbia MO 65211}
\email{edidin@math.missouri.edu}
\author{William Graham}
\address{Department of Mathematics\\ University of Georgia\\
Boyd Graduate Studies Research Center\\Athens, GA 30602}
\email{wag@math.uga.edu}
\thanks{Both authors were partially supported by the NSF}
\maketitle
\centerline{\it Dedicated to William Fulton on his 60th birthday}
\section{Introduction}

The purpose of this paper is to provide a characterization of
solvable linear algebraic groups in terms of a geometric property
of representations.  Representations with a related property
played an important role in the proof of the equivariant
Riemann-Roch theorem \cite{ERRoch}.  In that paper, we
constructed representations with that property (which we
call freely good) for the group of upper triangular matrices
in $GL_n$.  We noted that it seemed unlikely that such representations
exist for arbitrary groups; the main result of this paper implies
that they do not.

To state our results, we need some definitions.  A representation
$V$ of a linear algebraic group $G$ is said to be {\it good } (resp.
{\it freely good})
if there exists a non-empty
$G$-invariant open subset $U \subset V$ such that

(i) $G$ acts properly (resp. freely) on $U$.

(ii) $V \smallsetminus U$ is the union of a finite
number of $G$-invariant linear subspaces.

Note that freely good representations were called good
in \cite{ERRoch}.

The main result of the paper is the following theorem.

\begin{theorem} \label{A}
Let $G$ be a connected algebraic group over a field $k$ of
characteristic not equal to $2$.  Then $G$ is solvable if and only if
$G$ has a good representation.  Moreover, if $G$ is solvable and $k$
is perfect then $G$ has a freely good representation.
\end{theorem}

In characteristic $2$, a solvable group still has good representations,
and a partial converse holds (Corollary \ref{c.char2}).
A key step in the proof of the main result is Theorem \ref{t.nt}, which
is inspired by an example of Mumford \cite[Example 0.4]{GIT}.

In characteristic 0 solvable groups are characterized by a
weaker property which does not require the action to be proper.
(In general, if $G$ acts properly on $X$, then $G$ acts with
finite stabilizers on $X$, but the converse need not hold.)

\begin{theorem} \label{B}
Let $G$ be a connected algebraic group over
a field of characteristic $0$. Suppose that $G$ has a representation
$V$ that contains a nonempty open set $U$ such that

$(1)$ The complement of $U$ is a finite union of invariant linear subspaces,
and

$(2)$ $G$ acts with finite stabilizers on $U$.

\noindent Then $G$ is solvable.
\end{theorem}

Examples (see Section \ref{s.examples}) show that this weaker property
does not characterize solvability in positive characteristic.

\section{Preliminaries}

{\bf Groups and representations} We let $k$ denote a field, with
algebraic closure $\kbar$ and separable closure $k_s$. If $Z$ is a
$k$-variety and $k' \supset k$ is any extension of $k$ then $Z(k')$
denotes the $k'$-valued points of $Z$, while $Z_{k'}$ denotes the
$k'$-variety $Z \times_k k'$.

All groups in this paper are assumed to be linear algebraic groups
over a field $k$.  We assume that such a group $G$ is geometrically
reduced, that is, that $G_{\kbar}$ is reduced.  The identity component
of a group $G$ is denoted $G^0$.

Unless otherwise stated, a representation $V$ of a group $G$ is
assumed to be $k$-rational; i.e. $V$ is a $k$-vector space
and the action map $G \times V \to V$ is a morphism of $k$-varieties.

If $k' \supset k$ is a field extension we call
a $k'$-rational representation $V$
of $G_{k'}$ a
$k'$-representation of $G$.
We say that $V$ is defined
over $k$ if it is obtained by base change
from a $k$-rational representation.

If $k' \supset k$ be a Galois field extension, then $\Gal(k'/k)$ acts
on $k'$-representations of $G$. Indeed, let $V$ be a
$k'$-representation of $G$ corresponding to a $k'$ morphism $\rho
\colon G_{k'} \to GL(V)$.  For $g \in G(k_s)$ $^\sigma \rho(g)$ is
defined as follows (cf. \cite[AG14.3, 24.5]{Borel}).  Because $\rho$
is defined over $k'$, for any $\tau \in \Gal(k_s/k')$ and any $g \in
G(k_s)$, $\tau(\rho(\tau^{-1}(g))) = \rho(g)$.  Thus, if $\sigma \in
\Gal(k'/k) = \Gal(k_s/k)/\Gal(k'/k)$, then
$$\sigma ' (\rho((\sigma ')^{-1} g)) \in GL(V)(k_s)$$
is independent of the lift of $\sigma$ to an element
$\sigma ' \in \Gal(k_s/k)$. We will call this
point $\sigma \rho(\sigma^{-1}(g))$ and set
$\rho(g) = \sigma(\rho(\sigma^{-1} g))$.

The $k'$-representation $V$ is obtained by base change from
a representation defined over $k$ if and only if
$^\sigma \rho = \rho$ for all $\sigma \in \Gal(k'/k)$.

\medskip

{\bf Free and proper actions}
The action of a group
$G$ on a scheme $X$ is said to be {\it free}
if the action map $G \times X \to X \times X$ is
a closed embedding. The action is said to be
{\it proper} if the map $G \times X \to X \times X$
is proper.
If the action is proper then the stabilizer of
every point is finite. If the stabilizer
of every geometric point is a trivial group-scheme then we say
that the action is {\it set theoretically free}.
An action which is set theoretically free
and proper is free \cite{EIT}.

Let $H \to G$ be a finite morphism of algebraic groups.
If $G$ acts properly on a scheme $X$
then $H$ also acts properly on $X$. Thus, if $V$
is a good representation of $G$ then $V$ is also
a good representation of $H$ via the action induced
by the map $H \to G$. Moreover, if $H$ is a closed subgroup
and $V$ is a freely good representation of $G$, then
$V$ is a freely good representation of $H$.

\medskip

\begin{example}
Let $B$ be the group of upper triangular matrices in $\GL(n)$.  The
group $B$ acts by left multiplication on the vector space $V$ of
upper triangular matrices; it acts with trivial stabilizers on the
open subset $U$ of invertible upper triangular matrices.  Since the
matrices are upper triangular, $V \smallsetminus U$ is the union of
the invariant subspaces $L_i = \{ A \in V| A_{ii} = 0\}$.  This representation
is freely good because the action of $B$ on $U$ is identified
with $B$ acting itself by left multiplication. The map
$B \times B \to B \times B$ given by $(A,A') \mapsto (A,AA')$
is an isomorphism, so the action of $B$ on $U$ is free.

By contrast, the action
of $\GL(n)$ by left multiplication on the vector space $M_n$ of
$n \times n$ matrices is not good.
\end{example}

\section{Existence of good representations}
In this section we show that every connected solvable group has
good representations, and if $k$ is perfect,
freely good representations.

By the Lie-Kolchin
theorem $G_{\kbar}$ is trigonalizable; i.e.,
it can be embedded in the group $B_{\kbar} \subset \GL_n$
of upper triangular matrices.

Let $V_{\kbar}$ be the vector space of upper triangular $n \times n$
matrices.  The group $B_{\kbar}$ acts on $V_{\kbar}$ by left
multiplication and we have seen that this representation is freely
good. By restriction $V_{\kbar}$ is a good representation of
$G_{\kbar}$.  Consider the morphism $\rho\colon G_{\kbar} \rightarrow
\GL(V_{\kbar})$ corresponding to the action of $G_{\kbar}$ on
$V_{\kbar}$.

Since $\rho$ is a morphism of schemes of finite type, it is defined
over a field extension $k' \supset k$ of finite degree.  Write $V =
V_{k'}$ for the corresponding $k'$-representation; then we have $\rho:
G_{k'} \rightarrow \GL(V)$.

\medskip

{\bf Case I.} $k'$ is separable over $k$. (This will
occur when $k$ is perfect.)
In this case we will use Galois descent to construct
a freely good representation of $G$.

Replacing $k'$ by a possibly bigger field extension we may assume that
$k' \supset k$ is Galois.  Enumerate the elements of $\Gal(k'/k)$ as
$\{1 =\sigma_1, \sigma_2 , \ldots \sigma_d\}$ and consider the
representation $\Phi\colon G_{k'}\to \GL(V^{\oplus d})$ where $G_{k'}$
acts on the $j$-th factor by the representation $^{\sigma_j}\rho
\colon G_{k'} \to \GL(V)$.

We define $U_d \subset V^{\oplus d}$ to be the open set whose
$k_s$-rational points are the $d$-tuples
$(A_1, \ldots, A_d)$ where some $A_i$ is invertible.
We realize $U_d$ as a complement of $G_{k'}$-invariant
linear subspaces as follows. Let $L_j = \{A \in V| A_{jj} =  0\}$,
a $G_{k'}$-invariant subspace of $V$.  Given a $d$-tuple
$(j_1, \ldots, j_d)$, define
$$
L_{(j_1, \ldots, j_d)} = L_{j_1} \oplus \cdots \oplus
L_{j_d}.
$$ This is a $G_{k'}$-invariant subspace of $V^{\oplus d}$
and
$U_d = V_d \smallsetminus  \bigcup L_{(j_1, \ldots, j_d)}$.

\begin{lemma}
(cf. \cite[Theorem 2.2]{ERRoch}) $G_{k'}$ acts freely on $U_d$.
\end{lemma}
\begin{proof}
Since $G_{k'}$ is a closed subgroup of $B_{k'}$ and the open
set $U_d$ is $B_{k'}$ invariant, it suffices
to show that $B_{k'}$ acts freely on $U_d$.

To do this we must show that the map $B_{k'} \times U_d \to U_d \times U_d$
given on $k_s$ points by
$$(A,A_1, \ldots , A_d) \mapsto (AA_1, \sigma_2(A\sigma_2^{-1}(A_2)),
\ldots, \sigma_d(A \sigma_d^{-1}(A_d))). $$
is a closed embedding.

First we show that the image $Z$ of $B_{k'} \times U_d$ is
closed in $U_d \times U_d$.
Let $(A_1, A_2, \ldots A_d, C_1, \ldots , C_d)$ be matrix coordinates
on $U_d \times U_d$.
Expanding the inverse out in terms of the adjoint
we see that the image is contained in the subvariety defined
by the matrix equations
$$\sigma_j \sigma_i^{-1}(\det A_i) C_j = \sigma_j\sigma_i^{-1}(C_i \mbox{ }\Adj
A_i) A_j.$$
Suppose that a 2d-tuple of matrices
$(A_1, A_2, \ldots , A_d, C_1, C_2, \ldots C_d) \in U_d \times U_d$
satisfies the matrix equations above. At least one of
the $A_i$ and one of the $C_j$ is invertible because
we are in $U_d \times U_d$.
Let $A = \sigma_i^{-1}(C_iA_i^{-1})$. Substituting
into our equations we see that $C_l = \sigma_l(A\sigma_l^{-1}(A_l))$
for all $l$. Moreover $A$ is invertible since
$C_j$ is invertible and $C_j = \sigma_j(A \sigma_j^{-1}(A_j))$.
Hence every point satisfying the matrix equations is in
the image $Z$ of $B_{k'} \times U_d$, so $Z$ is closed.

The variety $Z$ is covered by open sets
of the form
$$
\{(A_1,A_2, \ldots A_j, \ldots A_d, AA_1, \ldots , \ldots
AA_j, \ldots A_d)| \det A_j \neq 0\}.
$$ These open sets
are isomorphic to $V^{d-1} \times B_{k'}$, where
$V^{d-1}$ is the $d-1$-fold cartesian product of
$V$. Hence the image
is smooth, in particular, normal.
The action of $G_{k'}$ on $U_d$ is set-theoretically free
so $G_{k'} \times U_d \to Z$ is a birational bijection.
By Zariski's main theorem (cf. \cite[AG18]{Borel})
a birational bijection of  a normal varieties is an isomorphism,
so $G_{k'} \times U_d \to Z$ is an isomorphism. Therefore,
$G_{k'} \times U_d \to U_d \times U_d$ is a closed embedding.
\end{proof}

\medskip

\textbf{Remark.} The proof of \cite[Theorem 2.2]{ERRoch} is
incomplete; the last paragraph of the above argument is
needed.

\medskip

For any basis of $V$, there is a natural choice of basis so that
with respect to this basis, if $g \in G(k_s)$, $\Phi(g)$
is represented
by the block diagonal matrix
$$
\begin{bmatrix}
\rho(g) & & & &\\ & ^{\sigma_2}\rho(g) & & & \\
& & .. & &\\ & & & .. &\\ & & & & ^{\sigma_d}\rho(g)
\end{bmatrix}$$
This representation is not defined over $k$ because the Galois group
acts by permuting the blocks.  More precisely, we have the following.
Given
a $d \times d$ matrix $M$, let $M[n]$ denote the $nd \times nd$ matrix
whose $ij$ block is $M_{ij} \cdot I_n$, where $I_n$ is the
$n \times n$ identity matrix.  If $\sigma \in \Gal(k'/k)$,
let $J_{\sigma}$ denote the permutation matrix corresponding
to the permutation $\sigma_i \mapsto \sigma \sigma_i$.
In matrix form,
for $g \in G(k_s)$,
$$
^\sigma \Phi(g) = J_{\sigma}[n]^{-1} \Phi(g) J_{\sigma}[n].
$$

We will show that $\Phi$ is $k'$-isomorphic to a freely good
representation defined over $k$.  Choose a primitive element $\alpha$
for the extension $k' \supset k$, and let $A$ be the $d \times d$
matrix with $A_{ij} = \sigma_j(\alpha^i)$.  $A$ is invertible since
$\alpha, \sigma_2(\alpha), \ldots \sigma_d(\alpha)$ are exactly the
roots of the irreducible polynomial $f \in k[x]$ of $\alpha$ over $k$,
so $\det A = \prod_{i < j} (-1)^{d}(\sigma_i(\alpha) -
\sigma_j(\alpha)) \neq 0$.  The Galois group acts by $\sigma(A) = A
J_{\sigma}$.  Consider the morphism $\Psi: G_{k'} \rightarrow
\GL(V^{\oplus d})$ defined by by $\Psi(g) = A[n] \Phi(g) A[n]^{-1}$
for $g \in k_s$ .  Then $^\sigma \Psi(g) = \Psi(g)$ for any $\sigma
\in \Gal(k'/k)$; hence $\Psi$ is defined over $k$.
Each of the
subspaces $L_{(j_1, \ldots, j_d)}$ is $G_{k'}$-invariant under the
action $\Psi$.  Moreover, because each
$L_{(j_1, \ldots , j_d)}$ is a vector subspace of
$V^{\oplus d}$ defined over $k$, the corresponding sub-representations $G_{k'}
\rightarrow \GL(L_{(j_1, \ldots, j_d)})$ are also defined over $k$.
Therefore $\Psi$ is obtained by base change from a freely good
representation of $G$.

\medskip

{\bf Case II. The general case}
In this case we may assume that there is a freely good $k'$-rational
representation
$\rho:G_{k'} \rightarrow \GL(V)$ defined
over a finite normal extension $k'$ of $k$.
Then $k' \supset k$ factors as $k' \supset k'' \supset k$, with
$k'/k''$ purely inseparable of degree $p^n$ and $k''/k$ Galois.
The Frobenius endomorphism on $V$ induces a group homomorphism
of $\GL(V)$. Composing $\rho$ with the $n$-th power of Frobenius
on $\GL(V)$ we obtain a representation defined over $k''$.
Because the Frobenius has finite kernel, this representation will
no longer be faithful.
However, the action of Frobenius is trivial on geometric points,
so $G$ will act properly on an open set whose
complement is a union of linear subspaces.
We can now use the Galois descent argument of Case I to obtain
a good $k$-rational representation of $G$.

\section{Characterization of solvable groups by good representations}
\label{s.reductive}

In this section we show that if $\mbox{char }k \neq 2$, every reductive
group with a good representation is a torus.  However,
many of the results of this section are valid in arbitrary
characteristic, and we only need that $\mbox{char }k \neq 2$ in part of
the proof of Theorem \ref{A}.  We will explicitly say when we
start assuming this; until then $\mbox{char }k$ is arbitrary.

Let $T$ be the diagonal torus in $SL_2$ and $N(T)$ the normalizer of $T$.
We will first show that $N(T)$ has no good representations.
We begin by recalling some facts about $N(T)$.
Let
$$
J =
\begin{bmatrix}
0 & -1 \\
1 & 0
\end{bmatrix}.
$$
We will also write
$$
H(t) =
\begin{bmatrix}
t & 0 \\
0 & t^{-1}
\end{bmatrix}.
$$
The group $N(T)$ is generated by $T$ and $J$; it has two components, $T$ and
$J(T)$.
The action of $SL_2$ on its two dimensional standard representation
$V$ induces an action on $S(V^*) \cong k[x,y]$, given by
$$
\begin{bmatrix}
a & b \\
c & d
\end{bmatrix}
: y \mapsto ay - cx , \qquad x \mapsto -by + dx.
$$
Let $W_i$ denote the subspace of $S(V^*)$ spanned by $x^i$ and $y^i$;
this is an irreducible representation of $N(T)$, of dimension $2$ (if
$i >0$).  Let $W'_0$ denote the $1$-dimensional irreducible
representation of $N(T)$ on which $T$ acts trivially and $J$ acts by
multiplication by $-1$.

If $\mbox{char }k \neq 2$,
the group $N(T)$ is linearly reductive;
that is, its action on any representation is completely reducible
\cite[p.191]{GIT}.  The next lemma shows that much of this
survives in arbitrary characteristic.

\begin{lemma}
Let $V$ be a representation of $N(T)$.

$(1)$ As a representation of $N(T)$, $V$ splits as a direct sum
of $N(T)$-submodules:
$$
V = V_0 \oplus \bigoplus_{i>0} V_{\pm i}.
$$
Here $V_{\pm i}$ is the sum of the $i$ and $-i$ weight spaces of $T$ on
$V$, and $V_j$ is the $j$-weight space.

$(2)$ The action of $N(T)$ on $V_{\pm i}$ $(i>0)$ is completely reducible,
and $V_{\pm i}$ is isomorphic as $N(T)$-module to a direct sum
of copies of $W_i$.

$(3)$ If $\mbox{char }k \neq 2$, then $V_0$ is isomorphic
to a direct sum of copies of $W_0$ and $W'_0$.
\end{lemma}

\begin{proof}
(1) Because the
action of $T$ on $V$ is completely reducible, we can decompose $V
= \oplus V_i$ as $T$-module.  As $J H(t) J^{-1} = H(t^{-1})$, we
have $J V_i = V_{-i}$.  Hence $V_{\pm i}$ is an $N(T)$-submodule
and we get the desired direct sum decomposition of $V$.

(2) Let $v_1, \ldots, v_d$ be a basis for $V_i$ $(i>0)$.
The map $v_r \mapsto y, J v_r
\mapsto -x$ defines an isomorphism of the span of
$v_r, Jv_r$ (denoted $\langle v_r, J
v_r \rangle$) with $W_i$, and the map
$W_i^{\oplus d} \rightarrow V_{\pm i}$, taking
the $r$-th component to $\langle v_r, J
v_r \rangle$, is an $N(T)$-module isomorphism.

(3) Decompose the $0$-weight space of $V$ into
the $+1$ and $-1$ eigenspaces of $J$; these are isomorphic to sums of
copies of $W_0$ and $W'_0$, respectively.
\end{proof}

The proof of the following result was motivated by
\cite[Example 0.4]{GIT}.

\begin{theorem} \label{t.nt}
The group  $N(T)$ has no
good representations.
\end{theorem}

\begin{proof}
If a group $G$ has good representations, then so does $G_{\bar{k}}$,
so we may assume that $k$ is algebraically closed.  Suppose that $V$
is a representation of $N(T)$, and let $U \subset V$ be the complement
of a finite set of invariant linear subspaces $\cS$.  We will show
that $N(T)$ does not act properly on $U = V - \cup_{L \in \cS}L$.  The
strategy of the proof is as follows.  Consider the action map $\Phi:
N(T) \times U \rightarrow U \times U$.  We will find a closed
subvariety $Z$ of $N(T) \times U$ whose closed points are of the form
$$
(
\begin{bmatrix}
0 & -\lambda^{-1} \\
\lambda & 0
\end{bmatrix}
, v_{\lambda})
$$
whose image is not closed in $U \times U$.
Hence $\Phi$ is not proper, so the representation is not good.

We now carry out the proof.  Decompose $V = \oplus V_i$,
where $V_i$ is the $i$-weight space of $V$ for $T$.  Pick $u \in U$, and write
$u = \sum u_i$ where $u_i \in V$.  Some of the $u_i$ may be $0$; let
$d$ be the dimension of the space spanned by the nonzero $u_i$.

\medskip

{\em Step 1.} If $a_i \neq 0$ for all $i$ with $u_i \neq 0$, then
$w = \sum a_i u_i \in U$.  Indeed, suppose not; then $w \in L$ for
some $L \in \cS$.  For almost all choices
$t_1, \ldots, t_d$ of $d$ elements of $k^*$, the vectors
$$
H(t_q) w = \sum_p t_q^{i_p} a_{i_p} u_{i_p} \in L
$$
are linearly independent.  (Here $i_1, \ldots, i_d$ are the indices
$i_p$ with $u_{i_p} \neq 0$.)  This follows because
the $d \times d$ matrix $A$ with entries
$$
A_{pq} = t_q^{i_p}
$$
is nonsingular for almost all $t_1, \ldots, t_d$.  (This
is because $\det A$ is a sum of monomials, where each monomial is a
product of one term from each row and each column.  Each monomial
has different multi-degree, so $\det A$ is not the zero polynomial.)
Therefore, the vectors $H(t_q) w$ span the same space as the $u_i$, so
$u \in L$, contradicting our assumption that $u \in U$.  We conclude
that $w \in U$, as claimed.

A similar argument shows that $J u_0 + \sum_{i \neq 0} u_i \in U$.
\medskip

{\em Step 2.} There exists an element $u' = \sum u'_i \in U$
with
$J u'_i = u'_{-i}$ for all $i>0$.  To see this, suppose
$u_{-j} \neq J u_j$ for some $j>0$.  Let $W_j \subset V_j \oplus V_{-j}$
be the subspace of vectors of the form $v_j + J v_j$ $(v_j \in V_j)$.
Note that $W_j$ generates $ V_j \oplus V_{-j}$ as $N(T)$-module.
Consider the affine linear subspace
$$
B = \sum_{i \neq \pm j} u_i + W_j.
$$
We claim that $B \cap U$ is nonempty.  If it is empty, then because
$B$ is affine linear and is contained in a finite union of the subspaces
in $\cS$, we see that $B \subset L$ for some $L \in \cS$.  But then
the span of $B$ is contained in $L$, so $\sum_{i \neq \pm j} u_i \in L$
and $W_j \subset L$.  Because $L$ is $N(T)$-stable,
$V_j \oplus V_{-j} \subset L$ as well; so
$$
u \in \sum_{i \neq \pm j} u_i + (V_j \oplus V_{-j}) \subset L,
$$
contradicting $u \in U$.  We conclude that $B \cap U$ is nonempty.
Replacing $u$ by an element of $B \cap U$, which we again
call $u$, we do not change $u_i$ for $i \neq \pm j$, but
we obtain $u_{-j} = J u_j$.  Iterating this process, we obtain $u'$
of the desired form.

Replacing $u$ by $u'$, we will assume that $J u_i = u_{-i}$ for all $i>0$.
 From the $N(T)$-module isomorphism of $\langle u_i, J u_i \rangle$
with $\langle x^i, y^i \rangle$, we see that for $i >0$,
$$
\begin{bmatrix}
0 & - \lambda^{-1} \\
\lambda & 0
\end{bmatrix}
: u_i \mapsto \lambda^i u_{-i}, \qquad u_{-i} \mapsto (-1)^i \lambda^{-i} u_i.
$$
Note also that
$$
\begin{bmatrix}
0 & - \lambda^{-1} \\
\lambda & 0
\end{bmatrix}
u_0 = J u_0.
$$

\medskip

{\em Step 3}.  Define
$$
v_{\lambda} = u_0 + \sum_{i>0}(u_{-i} + \lambda^{-i} u_i) \\
$$
$$
v'_{\lambda} = Ju_0 + \sum_{i>0}(u_{-i} + (-1)^i \lambda^{-i} u_i).
$$
For all $\lambda \neq 0$, both $v_{\lambda}$ and $v'_{\lambda}$
are in $U$ (by Step 1).  Define $Z$ to be the closed subvariety
of $N(T) \times U$ whose points are the pairs
$$
(
\begin{bmatrix}
0 & - \lambda^{-1} \\
\lambda & 0
\end{bmatrix}
, v_{\lambda}).
$$
Then
$$
\Phi(Z) = \{ (v_{\lambda}, v'_{\lambda}) \}_{\lambda \neq 0} \subset U \times
U.
$$
Consider the point
$$
(v, v') = ( u_0 + \sum_{i>0} u_{-i}, Ju_0 + \sum_{i>0} u_{-i}).
$$
Reasoning as in
Step 1 shows that $u$ is in the $N(T)$-module generated by $v$ or $v'$,
so if either $v$ or $v'$ were in $L$ then $u$ would be, but this
is impossible as $u \in U$.  Hence $v$ and $v'$ are in $U$,
so $(v, v') \in U \times U$.  Also, $(v, v')$ is not in
$\Phi(Z)$, but is in the closure of $\Phi(Z)$ in $U \times
U$.  We conclude that $\Phi$ is not proper, so the
representation is not good.
\end{proof}

\subsection{Proof of Theorem \ref{A}}
Let $G$ be a connected nonsolvable linear algebraic group.  Consider the
surjective map $\pi: G \rightarrow G_1 = G/\R_uG$, where $\R_uG$ is
the unipotent radical of $G$ and $G_1$ is reductive.  Because $G$ is not
solvable, $G_1$ is not
trivial or a torus.  Let $T$ be a
maximal torus of $G$.  Then $T_1 = \pi(T)$ is a maximal torus of
$G_1$, and $\pi$ induces an isomorphism of Weyl groups $W(T, G)
\rightarrow W(T_1,G_1)$ \cite[11.20]{Borel}.
(Here $W(T, G) = N_G(T)/Z_G(T)$ where $N_G$ and $Z_G$ denote normalizer and
centralizer of $T$ in $G$, and similarly for $G_1$.)
Because $\ker \pi$ is a
unipotent group, $\pi|_T: T \rightarrow T_1$ is an isomorphism.
Note that $Z_G(T) = T
\cdot (\R_uG)^T$ \cite[13.17]{Borel}.  Because $G_1$ is reductive, this
fact (applied to $G_1$) implies that $Z_{G_1}(T_1) = T_1$.
Moreover, any $g_1 \in N_{G_1}(T_1)$ can be lifted to
$g \in N_G(T)$.  This follows because the isomorphism of Weyl groups
above, and the structure of the centralizers, imply that each
component of $N_{G_1}(T_1)$ is the image of a surjective map of a
component of $N_G(T)$.

As $G_1$ is not a torus, there is a root $\alpha$ and a homomorphism
$\phi_{\alpha}: \SL_2 \rightarrow G_1$ with kernel either trivial, or
the set of matrices
$\begin{bmatrix} a & 0 \\ 0 & a \end{bmatrix}$ with $a^2 = 1$.
Moreover (using the subscript $\SL_2$ to denote terms for $\SL_2$
defined in the previous subsection), $\phi_{\alpha}(T_{\SL_2}) \subset
T$ and $J_1 := \phi_{\alpha}(J_{\SL_2}) \in N_{G_1}(T_1)$.  (See
\cite[p.176]{Jantzen} for these facts.)  Let $H_1 =
\phi_{\alpha}(N(T_{\SL_2}))$; its identity component $H_1^0 =
\phi_{\alpha}(T_{\SL_2}) \subset T_1$.  Because $H_1$ is a finite image
of $N(T_{\SL_2})$, it has no good representations (and hence neither
does $G_1$).

{\it Up to this point, $\mbox{char }k$ has been arbitrary; now we
assume that $\mbox{char }k \neq 2$.}

Because $\pi|_T$ is an isomorphism
there is a unique subgroup $H^0\subset T$ projecting isomorphically
to $H_1^0$.  As noted above, we can choose a lift $J \in N_G(T)$
of $J_1 \in N_{G_1}(T_1)$.  Write $J = J_s J_u$ for the Jordan
decomposition of $J$.  Because $\mbox{char }k \neq 2$,
$J_1$ is semisimple, so $\pi(J_u) = 1$.  Therefore we
can replace $J$ by $J_s$ and assume $J$ is semisimple.
Now, $J^2$ corresponds to the identity element
in the Weyl group (as $J_1^2$ does), so $J^2 \in
Z_G(T) = T \cdot (\R_uG)^T$. Since $J$ is semisimple, we
conclude that $J^2 \in T$.  As $J_1^2$ is in the subgroup
$H_1^0$ of $T_1$ and $T$ maps isomorphically to $T_1$, we
conclude that $J^2 \in H^0$.  Therefore the group
$H$ generated by $H^0$ and $J$ maps isomorphically to $H_1$, and thus
has no good representations.  Therefore
$G$ has no good representations.  This proves Theorem \ref{A}.

\medskip

The proof of Theorem \ref{A} yields the following weaker
statement in characteristic $2$.  Note that Levi decompositions
need not exist in positive characteristic \cite[11.22]{Borel}.

\begin{corthm} \label{c.char2}
Suppose  $\mbox{char } k = 2$.  If the connected algebraic group
$G$ has a Levi decomposition
and if $G$ has a good representation, then $G$ is
solvable. In particular, any connected reductive group with a good
representation is diagonalizable.
\end{corthm}

\begin{proof}
Suppose $G = LN$ where $L$ is reductive and $N$ unipotent.  If
$G$ has a good representation then so does $L$.  As proved above,
this implies that $L$ is a torus, so $G$ is solvable.
\end{proof}

\section{Proof of Theorem \ref{B}}
If $G$ has a representation
$V$ that contains an open
set $U$ whose complement is a finite union of invariant subspaces
such that $G$ acts with finite stabilizers on $U$, then
$G_{\kbar}$ also has such a representation. Thus we can
assume that $k$ is algebraically closed.

Assume that $G$ is not solvable and let $V$ be a representation of
$G$. Since the characteristic is 0, $G$ has a Levi subgroup $L$.
Since $G$ is assumed to be non-solvable, $L$ contains a Borel
subgroup which is not a torus. Hence $L$
contains a non-trivial unipotent subgroup $N$.

Since the characteristic is 0 and $L$ is reductive, $V$ decomposes as
a direct sum $V = V_1 \oplus V_2 \ldots V_p$ of irreducible
$L$-modules.  Every vector in the subspace $V^N = V_1^N \oplus V_2^N
\oplus \ldots \oplus V_p^N$ has positive
dimensional stabilizer.  Since $N$ is unipotent, $V_i^N \neq 0$ for
each $i$, so $L(V_i^N)$ spans all of $V_i$.
Hence the subset $LV^N = L(V_1^N \oplus \ldots \oplus
V_p^N)$ which consists of vectors with positive dimensional
stabilizers cannot be contained in any proper $L$-invariant
subspace. Since $L$ is a subgroup of $G$ this means that $LV^N$ is not
contained in any proper $G$-invariant subspace.
Hence $V$ does not have properties (1) and (2). \endproof

\section{Examples and complements} \label{s.examples}
In this section we discuss ``set-theoretic'' versions
of the conditions freely good and good.
We will say a representation $V$ is set-theoretically
freely good (resp. set-theoretically good)
if it contains a nonempty open subset $U$ whose complement
is a union of invariant subspaces, such that $G$ acts
with trivial stabilizers (resp. finite stabilizers)
on $U$ (cf. Theorem \ref{B}).  Surprisingly, these conditions are not
enough to characterize solvability in arbitrary characteristic.

\begin{example}
Let $V$ be the standard representation
of $\SL_2$ and let $V_d = S(V^*)$ be the vector space of homogeneous forms
of degree $d$. As in Section \ref{s.reductive},
$\SL_2$ acts on $V_d$.  If $p = \mbox{char }k$ is an odd prime
then $W_p = V_{2p-2} \oplus V_1$ is a set-theoretically freely
good representation of $\SL_2$. The reason is as follows.
The stabilizer of any pair
of forms $(f(x,y),l(x,y))$ is trivial as long as
$l(x,y) \neq 0$ and the coefficient of $x^{p-1}y^{p-1}$ in $f$ is
non-zero.
Since the characteristic is $p$,
the subspace $L_{2p-2} \subset V_{2p-2}$ of forms with no
$x^{p-1}y^{p-1}$ term is an $\SL_2$ invariant subspace
(cf. \cite[II2.16]{Jantzen}). Thus, $\SL_2$
acts with trivial stabilizers on the open
set $U_p  = W_p \smallsetminus (W_{2p-2} \oplus V_1 \cup V_{2p-2} \oplus 0)$.

In characteristic 2, the representation $W_2 = V_2 \oplus V_1$ is not
set-theoretically freely good because the matrix $\begin{bmatrix} 0 &
1\\1 & 0\end{bmatrix}$ stabilizes the pair $(x^2y^2, x+y)$.  However,
$W_2$ is set-theoretically good.

In positive characteristic,
we do not know if the group $\SL_n$ admits set-theoretically good
representations for $n \geq 3$.
\end{example}

\begin{example} Assume $k$ is algebraically closed and $\mbox{char }k \neq 2$.
Then $G = PGL_2$
has no representation which is set-theoretically freely good.  Indeed,
let
$$
g =
\begin{bmatrix}
0 & 1 \\
1 & 0
\end{bmatrix}
$$
and let $H = \{1, g\}$.  If $V$ is
any representation of $G$, $V^H$ generates $V$ as a representation
of $G$.  Indeed, this holds if $V$ is irreducible, since
for any vector $v$, the vector $v + gv$ is a nonzero $H$-invariant.
Because $\mbox{char }k \neq 2$, the action of $H$ is
completely reducible, so if
$$
0 \rightarrow V_1 \rightarrow V_2 \rightarrow V_3 \rightarrow 0
$$
is an exact sequence of $G$-modules, then the corresponding
sequence of $H$-invariants is also exact.  By induction,
we may assume that $V_1^H$ and $V_3^H$ generate $V_1$ and
$V_3$ as $G$-modules, and a diagram-chase then shows that
$V_2^H$ generates $V_2$ as a $G$-module.

Hence if $V$ is any representation, there is no proper
invariant linear subspace of $V$ containing $V^H$.
Therefore $V$ is not set-theoretically freely good.

A similar argument shows that $\PGL_n$ and $\GL_n$
do not have set-theoretically freely good representations.
\end{example}

We conclude with a proposition about the inductive construction
of good representations.

\begin{prop} \label{p.normal}
Let $G$ be a connected linear algebraic group and $H$ a normal
subgroup.  Assume $k$ is algebraically closed.  If $H$ and $G/H$ have
set-theoretically freely good representations, then so does $G$.

\end{prop}

\begin{proof}
For this proof only, we will use ``good'' to mean ``set-theoretically
freely good''.  Let $W$ be a good representation of $H$, with $M_i$ a
finite set of proper invariant subspaces containing the vectors with
nontrivial stabilizers.  Because $G/H$ is affine
\cite[Theorem 6.8]{Borel}, the vector bundle $G
\times^H W$ is generated by a finite dimensional space of global
sections $\Gamma$.  We will view sections of the vector bundle as
regular functions $\g: G \rightarrow W$ satisfying $\g(gh) = h^{-1}
\cdot \g(g)$, where on the right side we are using the action of $H$
on $W$.  The action of $G$ on the space of sections of the vector
bundle corresponds to the left action of $G$ on regular functions: $(g
\cdot \g)(g_0) = \g(g^{-1} g_0)$.  Because the action of $G$ on
regular functions is locally finite, by enlarging the space $\Gamma$,
we may assume $\Gamma$ is stable under the $G$-action.

Define $L_i$ to be the subspace of $\Gamma$ consisting of those elements
of $\Gamma$ which are sections of $G \times^H M_i$.  Each $L_i$ is
a $G$-stable subspace of $\Gamma$.  Let $\Gamma^0$ denote the
complement of the $L_i$ in $\Gamma$.

Let $V$ be a good representation of $G/H$, viewed as a representation
of $G$ via the map $G \rightarrow G/H$.  We claim that $V \oplus
\Gamma$ is a set-theoretically good representation of $G$.  Indeed,
let $V_j$ be a finite set of invariant subspaces of $V$ containing the
vectors with nontrivial stabilizer.  It suffices to show that the
vectors with nontrivial stabilizer in $V \oplus \Gamma$ are contained
in the union of the subspaces $V_j \oplus \Gamma$ and $V \oplus L_i$.
To see this, let $(v, \g)$ be in the complement of these subspaces. so
$v \notin V_j$ and $\g \notin L_i$ for any $i,j$.  We must show that
$\stab_G(v, \g)$ is trivial.  First, $\stab_G(v, \g) \subset
\stab_G(v) = H$.  Let $h \in \stab_G(v, \g)$.  As above, we will view
$\g$ as a function $G \rightarrow W$.  Because $\g$ is not in any
$L_i$, we have $\g$ is not a section of $G \times^H M_i$ for any $i$.
In other words, the open subsets $\g^{-1}(W \setminus M_i)$ of $G$ are
nonempty.  Choose $g_0$ in the intersection of these sets, so $s(g_0)
\notin M_i$ for any $i$.  Our hypothesis implies that $h \cdot \g =
\g$.  By definition, we have
$$
(h \cdot \g)(g_0) = \g(h^{-1}g_0) = \g(g_0 (g_0^{-1} h^{-1} g_0))
= (g_0^{-1} h g_0) \g(g_0).
$$
But $\stab_H \g(g_0) = \{ 1 \}$, so
we conclude $h = 1$ as desired.
\end{proof}


\begin{thebibliography}{99}
\bibitem[Borel]{Borel} A. Borel, {\it Linear Algebraic Groups} (2nd enlarged
edition),
G.T.M {\bf 126}, Springer-Verlag (1991).

\bibitem[EG1]{EIT} D. Edidin, W. Graham, {\it Equivariant intersection
theory}, Invent. Math., {\bf 131} (1998), 595-634.

\bibitem[EG2]{ERRoch} D. Edidin, W. Graham, {\it Riemann-Roch for equivariant
Chow groups}, Duke Math. J., to appear.

\bibitem[GIT]{GIT} D. Mumford, J. Fogarty, F. Kirwan, {\it Geometric
Invariant Theory: 3rd enlarged edition}, Springer-Verlag 1994.

\bibitem[J]{Jantzen} J. C. Jantzen, {\it Representations of Algebraic
Groups}, Academic Press (1987).

\end{thebibliography}
\end{document}